\documentclass{article}
\usepackage{amssymb}
\usepackage{amsfonts}
\usepackage{amsmath}

\input{tcilatex}
\begin{document}

\begin{flushright}
{\huge Eulerian Composition of Certain Franklin Squares}{\LARGE \bigskip }

{\Large RONALD P. NORDGREN}

{\small Rice University\vspace{-0.03in}}

{\small Houston, TX 77251\vspace{-0.03in}}

{\small nordgren@rice.edu}
\end{flushright}

\begin{quote}
\noindent \textbf{Abstract.}{\small \ Several specific Franklin squares and
magic squares are decomposed into their quotient and remainder squares. The
results support the conjecture that Franklin used the Eulerian composition
method to construct many of his squares. This method also can be used to
construct new Franklin squares as illustrated herein.}\smallskip 
\end{quote}

\noindent {\Large Introduction}\smallskip

Franklin squares were created by Benjamin Franklin and two of them were
first published in 1769 according to Pasles \cite[p. 127]{PASL}. They differ
from ordinary magic squares by the imposition of bent-diagonal sum
conditions instead of the two diagonal sum conditions. Franklin squares have
received considerable attention over the years as reviewed by Pasles \cite%
{PASL}. A general treatment of magic squares and circles is given by
Pickover \cite{PICK}.

Unfortunately, Franklin did not reveal his construction methods and there is
controversy on this question. His order-8 and order-16 squares of 1769, as
given in \cite{PASL}, can be constructed by a direct method due to Jacobs 
\cite{JACO} which he applies to construct an order-24 Franklin square. This
construction involves sequential placement of integer elements in five
prescribed steps. Jacobs' construction also yields the order-32 Franklin
square of Behrforooz \cite{BEHR}. However, a proof of the validity of
Jacobs' construction for order-$8k$ Franklin squares is lacking.

Franklin squares also can be constructed by the addition of two orthogonal
auxiliary squares according to a composition formula first published by
Euler \cite{EULER} in 1776. According to Pasles \cite{PASL}, the composition
construction was used by Youle in 1813 to construct Franklin's original two
squares and perhaps by Franklin himself (before Euler). However, it is not
clear how Franklin formed his two auxiliary squares, if indeed he used this
method.

Nordgren \cite{NORD1} gives a systematic construction for Franklin squares
of order $8k$. He obtains two types of formulas for the elements of the two
auxiliary squares that allow straightforward numerical formation of order-$%
8k $ Franklin squares. These formulas are used to verify that the auxiliary
squares are suitable for construction of Franklin squares. Since this
construction leads to Franklin's original 1769 squares for $k=1,2,$ it may
be regarded as an extension of his methodology and may have been used by
him. This speculation is supported by an examination of the Eulerian
composition of three order-8 Franklin squares by Pasles \cite{PASL} and
three order-16 Franklin squares in what follows. The auxiliary squares for
these Franklin squares follow regular patterns. The Eulerian composition of
some of these squares and others is given by Morris \cite{MORR} who also
contends that Franklin used this method. It is difficult to see how many of
these Franklin squares could be constructed by direct methods.

The formulas given in \cite{NORD1} also generate the order-24 Franklin
square of Jacobs \cite{JACO} (Figures 1,2,3) and the order-32 Franklin
square of Behrforooz \cite{BEHR}. A new order-40 Franklin square,
constructed by Nordgren's method, is shown in Figure 4. It also can be
constructed by Jacob's method. These constructions lend support to Jacobs'
construction but a direct connection between the two methods is lacking.

Numerical generation of all natural order-8 Franklin squares is carried out
by Schindel, et.al. \negthinspace \cite{SCH}. The basic 4320 of these
squares were made available to the author and one of them is decomposed into
auxiliary squares in what follows. \medskip

\noindent {\Large Definitions}\smallskip

For completeness, we review the definitions of magic and Franklin squares.
All rows and columns of a \emph{semi-magic} square matrix must sum to an 
\emph{index number} $m.$ If its main diagonal and the cross diagonal also
sum to $m,$ then the square is \emph{magic}. \emph{Natural} $n\times n$ (%
\emph{order}-$n)$ magic and semi-magic squares have elements $1,2,\ldots
,n^{2}$ for which%
\begin{equation}
m=\frac{n}{2}\left( n^{2}+1\right) \,.  \label{m}
\end{equation}%
Instead of the two diagonal sum conditions, Franklin defined four families
of \emph{bent diagonals} which must sum to $m$. In the following families of
order-6 squares, elements on the six bent diagonals (with wraparound) have
the same symbol:%
\begin{equation}
\begin{tabular}{|cccccc|}
\hline
$\!\spadesuit \!\!$ & $\!\!{\mathbf{\Diamond }}\!\!$ & $\!\!\clubsuit \!\!$
& $\!\!\heartsuit \!\!$ & $\!\!\blacktriangle \!\!$ & $\!\!\triangledown \!$
\\ 
$\!\triangledown \!\!$ & $\!\!\spadesuit \!\!$ & $\!\!{\mathbf{\Diamond }}%
\!\!$ & $\!\!\clubsuit \!\!$ & $\!\!\heartsuit \!\!$ & $\!\!\blacktriangle
\! $ \\ 
$\!\blacktriangle \!\!$ & $\!\!\triangledown \!\!$ & $\!\!\spadesuit \!\!$ & 
$\!\!\Diamond \!\!$ & $\!\!\clubsuit \!\!$ & $\!\!\heartsuit \!$ \\ 
$\!\blacktriangle \!\!$ & $\!\!\triangledown \!\!$ & $\!\!\spadesuit \!\!$ & 
$\!\!\Diamond \!\!$ & $\!\!\clubsuit \!\!$ & $\!\!\heartsuit \!$ \\ 
$\!\triangledown \!\!$ & $\!\!\spadesuit \!\!$ & $\!\!\Diamond \!\!$ & $%
\!\!\clubsuit \!\!$ & $\!\!\heartsuit \!\!$ & $\!\!\blacktriangle \!$ \\ 
$\!\spadesuit \!\!$ & $\!\!\Diamond \!\!$ & $\!\!\clubsuit \!\!$ & $%
\!\!\heartsuit \!\!$ & $\!\!\blacktriangle \!\!$ & $\!\!\triangledown \!$ \\ 
\hline
\end{tabular}%
\ \ \ 
\begin{tabular}{|cccccc|}
\hline
$\!\spadesuit \!\!$ & $\!\!{\mathbf{\Diamond }}\!\!$ & $\!\!\clubsuit \!\!$
& $\!\!\heartsuit \!\!$ & $\!\!\blacktriangle \!\!$ & $\!\!\triangledown \!$
\\ 
$\!{\mathbf{\Diamond }}\!\!$ & $\!\!\clubsuit \!\!$ & $\!\!\heartsuit \!\!$
& $\!\!\blacktriangle \!\!$ & $\!\!\triangledown \!\!$ & $\!\!\spadesuit \!$
\\ 
$\!\clubsuit \!\!$ & $\!\!\heartsuit \!\!$ & $\!\!\blacktriangle \!\!$ & $%
\!\!\triangledown \!\!$ & $\!\!\spadesuit \!\!$ & $\!\!{\mathbf{\Diamond }}%
\! $ \\ 
$\!\clubsuit \!\!$ & $\!\!\heartsuit \!\!$ & $\!\!\blacktriangle \!\!$ & $%
\!\!\triangledown \!\!$ & $\!\!\spadesuit \!\!$ & $\!\!{\mathbf{\Diamond }}%
\! $ \\ 
$\!\Diamond \!\!$ & $\!\!\clubsuit \!\!$ & $\!\!\heartsuit \!\!$ & $%
\!\!\blacktriangle \!\!$ & $\!\!\triangledown \!\!\!$ & $\!\!\spadesuit \!$
\\ 
$\!\spadesuit \!\!$ & $\!\!\Diamond \!\!$ & $\!\!\clubsuit \!\!$ & $%
\!\!\heartsuit \!\!$ & $\!\!\blacktriangle \!\!$ & $\!\!\triangledown \!$ \\ 
\hline
\end{tabular}%
\ \ \ 
\begin{tabular}{|cccccc|}
\hline
$\!\spadesuit \!\!$ & $\!\!\heartsuit \!\!$ & $\!\!\clubsuit \!\!$ & $%
\!\!\clubsuit \!\!$ & $\!\!\heartsuit \!\!$ & $\!\!\spadesuit \!$ \\ 
$\!\heartsuit \!\!$ & $\!\!\clubsuit \!\!$ & $\!\!\Diamond \!\!$ & $%
\!\!\Diamond \!\!$ & $\!\!\clubsuit \!\!$ & $\!\!\heartsuit \!$ \\ 
$\!\clubsuit \!\!$ & $\!\!\Diamond \!\!$ & $\!\!\blacktriangle \!\!\!$ & $%
\!\!\blacktriangle \!\!$ & $\!\!\Diamond \!\!$ & $\!\!\clubsuit \!$ \\ 
$\!{\mathbf{\Diamond }}\!\!$ & $\!\!\blacktriangle \!\!$ & $%
\!\!\triangledown \!\!$ & $\!\!\triangledown \!\!$ & $\!\!\blacktriangle
\!\! $ & $\!\!\Diamond \!$ \\ 
$\!\blacktriangle \!\!\!$ & $\!\!\triangledown \!\!$ & $\!\!\spadesuit \!\!$
& $\!\!\spadesuit \!\!$ & $\!\!\triangledown \!\!$ & $\!\!\blacktriangle \!$
\\ 
$\!\triangledown \!\!$ & $\!\!\spadesuit \!\!$ & $\!\!\heartsuit \!\!$ & $%
\!\!\heartsuit \!\!$ & $\!\!\spadesuit \!\!$ & $\!\!\triangledown \!$ \\ 
\hline
\end{tabular}%
\ \ \ 
\begin{tabular}{|cccccc|}
\hline
$\!\spadesuit \!\!$ & $\!\!\triangledown \!\!$ & $\!\!\blacktriangle \!\!$ & 
$\!\!\blacktriangle \!\!$ & $\!\!\triangledown \!\!$ & $\!\!\spadesuit \!$
\\ 
$\!\Diamond \!\!$ & $\!\!\spadesuit \!\!$ & $\!\!\triangledown \!\!$ & $%
\!\!\triangledown \!\!$ & $\!\!\spadesuit \!\!\!$ & $\!\!\!\Diamond \!$ \\ 
$\!\clubsuit \!\!$ & $\!\!\Diamond \!\!$ & $\!\!\spadesuit \!\!$ & $%
\!\!\spadesuit \!\!$ & $\!\!\Diamond \!\!$ & $\!\!\clubsuit \!\!$ \\ 
$\!\heartsuit \!\!$ & $\!\!\clubsuit \!\!$ & $\!\!\Diamond \!\!$ & $%
\!\!\Diamond \!\!$ & $\!\!\clubsuit \!\!$ & $\!\!\heartsuit \!$ \\ 
$\!\blacktriangle \!\!\!$ & $\!\!\heartsuit \!\!$ & $\!\!\clubsuit \!\!$ & $%
\!\!\clubsuit \!\!$ & $\!\!\heartsuit \!\!$ & $\!\!\blacktriangle \!$ \\ 
$\!\triangledown \!\!$ & $\!\!\blacktriangle \!\!$ & $\!\!\heartsuit \!\!$ & 
$\!\!\heartsuit \!\!$ & $\!\!\blacktriangle \!\!$ & $\!\!\triangledown \!$
\\ \hline
\end{tabular}%
\,.  \label{Bent}
\end{equation}%
Also, Franklin required that the elements on all left and right half-rows
and all upper and lower half-columns of his squares must sum to $m/2$ which
makes his squares semi-magic. This requirement forces natural Franklin
squares to be of doubly-even order $\left( n=4k\right) .$ In addition, he
required that elements in all $2\times 2$ subsquares (including broken ones)
of his order-$n$ square sum to $4m/n.$ Matrix formulas for the three
Franklin square conditions are given by Nordgren \cite{NORD2}. Hurkens \cite%
{HURK} shows by exhaustive search that no natural order-12 Franklin squares
exist. Also, Pasles \cite{PASL2} shows that there are no natural order-4
Franklin squares.

In a \emph{pandiagonal square}, the elements on all diagonals in both
directions (including broken ones) sum to $m.$ Some Franklin squares also
are pandiagonal, including one-third of the order-8 natural ones \cite{SCH}.
Nordgren \cite{NORD1} shows that a Franklin square of order $8k$ can be
transformed to a pandiagonal magic square in two ways but the converse is
not true in general. Also, Nordgren \cite{NORD2} shows that an order-$8k$
pandiagonal Franklin magic square can be obtained from transformation of a
complete (or most-perfect) magic square which are studied and enumerated by
Ollerenshaw and Br\'{e}e \cite{OLLE}.

For a square matrix $M_{n}$ of order-$n,$ Euler's composition formula \cite%
{EULER} can be written as%
\begin{equation}
M_{n}=nQ_{n}+R_{n}+U_{n}\,,  \label{Comb}
\end{equation}%
where $Q_{n}$ and $R_{n}$ are order-$n$ orthogonal matrices and $U_{n}$ is
the \emph{unity matrix} with all elements $1.$ \emph{Orthogonal matrices}
are defined as having each ordered pair of elements in the same position in
the two matrices occurring once and only once. A natural Franklin matrix can
be constructed by requiring that the \emph{quotient matrix} $Q_{n}$ and the 
\emph{remainder matrix} $R_{n}$ have elements $0,1,\ldots ,n-1$ repeated $n$
times. If such $Q_{n}$ and $R_{n}$ are orthogonal and satisfy the three
Franklin sum conditions with $m$ replaced by $\hat{m}=n\left( n-1\right) /2,$
then $M_{n}$ from (\ref{Comb}) is a Franklin matrix.

Given $M_{n},$ then $R_{n}$ and $Q_{n}$ are obtained from 
\begin{equation}
R_{n}=\left( M_{n}-U_{n}\right) \ \func{mod}\ n,\quad Q_{n}=\frac{1}{n}%
\left( M_{n}-U_{n}-R_{n}\right)  \label{RnQn}
\end{equation}%
which are used in what follows.\medskip

\noindent {\Large Order-6 Squares}\smallskip

Pasles \cite[p. 196]{PASL} gives the following square published by Franklin
in 1769:%
\begin{equation}
M_{6}=\left[ 
\begin{array}{cccccc}
2 & 9 & 4 & 29 & 36 & 31\vspace{-0.01in} \\ 
34 & 32 & 30 & 7 & 5 & 3\vspace{-0.01in} \\ 
6 & 1 & 8 & 33 & 28 & 35\vspace{-0.01in} \\ 
20 & 27 & 22 & 11 & 18 & 13\vspace{-0.01in} \\ 
25 & 23 & 21 & 16 & 14 & 12\vspace{-0.01in} \\ 
24 & 19 & 26 & 15 & 10 & 17\vspace{-0.01in}%
\end{array}%
\right] .  \label{M6}
\end{equation}%
It is natural, semi-magic, and its four main bent-diagonals (starting at the
corners) sum to $m=111$ as do the four bent diagonals that bend at the
middle of the four sides. In addition, other sets of squares add to $m$ \cite%
[p. 197]{PASL}. However, it does not satisfy all the Franklin square sum
conditions. Its auxiliary matrices from (\ref{RnQn}) are%
\begin{equation}
Q_{6}=\left[ 
\begin{array}{cccccc}
0 & 1 & 0 & 4 & 5 & 5\vspace{-0.01in} \\ 
5 & 5 & 4 & 1 & 0 & 0\vspace{-0.01in} \\ 
0 & 0 & 1 & 5 & 4 & 5\vspace{-0.01in} \\ 
3 & 4 & 3 & 1 & 2 & 2\vspace{-0.01in} \\ 
4 & 3 & 3 & 2 & 2 & 1\vspace{-0.01in} \\ 
3 & 3 & 4 & 2 & 1 & 2\vspace{-0.01in}%
\end{array}%
\right] ,\quad R_{6}=\left[ 
\begin{array}{cccccc}
1 & 2 & 3 & 4 & 5 & 0\vspace{-0.01in} \\ 
3 & 1 & 5 & 0 & 4 & 2\vspace{-0.01in} \\ 
5 & 0 & 1 & 2 & 3 & 4\vspace{-0.01in} \\ 
1 & 2 & 3 & 4 & 5 & 0\vspace{-0.01in} \\ 
0 & 4 & 2 & 3 & 1 & 5\vspace{-0.01in} \\ 
5 & 0 & 1 & 2 & 3 & 4\vspace{-0.01in}%
\end{array}%
\right] .  \label{QR6}
\end{equation}%
These auxiliary matrices are not semi-magic since their second and fifth
columns do not sum to $\hat{m}=15$ as do all other columns and all rows.
This counterexample indicates that semi-magic conditions on $Q_{n}$ and $%
R_{n}$ are not necessary conditions for a square to be semi-magic. Also,
some of the bent diagonals of $Q_{n}$ and $R_{n}$ mentioned above do not sum
to $\hat{m}=15$. However, a counterexample of a full Franklin square has not
been found and it is an open question whether Franklin conditions on $Q_{n}$
and $R_{n}$ are necessary for $M_{n}$ from (\ref{Comb}) to be a Franklin
square.

Here is the order-6 natural magic square given by Euler \cite{EULER}:%
\begin{equation}
M_{6}=\left[ 
\begin{array}{cccccc}
3 & 36 & 30 & 4 & 11 & 27\vspace{-0.01in} \\ 
22 & 13 & 35 & 12 & 14 & 15\vspace{-0.01in} \\ 
16 & 18 & 8 & 31 & 17 & 21\vspace{-0.01in} \\ 
28 & 20 & 6 & 29 & 19 & 9\vspace{-0.01in} \\ 
32 & 23 & 25 & 2 & 24 & 5\vspace{-0.01in} \\ 
10 & 1 & 7 & 33 & 26 & 34\vspace{-0.01in}%
\end{array}%
\right]   \label{E6}
\end{equation}%
with auxiliary matrices from (\ref{RnQn}):%
\begin{equation}
Q_{6}=\left[ 
\begin{array}{cccccc}
0 & 5 & 4 & 0 & 1 & 4\vspace{-0.01in} \\ 
3 & 2 & 5 & 1 & 2 & 2\vspace{-0.01in} \\ 
2 & 2 & 1 & 5 & 2 & 3\vspace{-0.01in} \\ 
4 & 3 & 0 & 4 & 3 & 1\vspace{-0.01in} \\ 
5 & 3 & 4 & 0 & 3 & 0\vspace{-0.01in} \\ 
1 & 0 & 1 & 5 & 4 & 5\vspace{-0.01in}%
\end{array}%
\right] ,\quad R_{6}=\left[ 
\begin{array}{cccccc}
2 & 5 & 5 & 3 & 4 & 2\vspace{-0.01in} \\ 
3 & 0 & 4 & 5 & 1 & 2\vspace{-0.01in} \\ 
3 & 5 & 1 & 0 & 4 & 2\vspace{-0.01in} \\ 
3 & 1 & 5 & 4 & 0 & 2\vspace{-0.01in} \\ 
1 & 4 & 0 & 1 & 5 & 4\vspace{-0.01in} \\ 
3 & 0 & 0 & 2 & 1 & 3\vspace{-0.01in}%
\end{array}%
\right] .  \label{EQR6}
\end{equation}%
This time $Q_{6}$ and $R_{6}$ are magic with $\hat{m}=15$.

Here is the order-6 natural magic square from the Historical Museum in Xian,
China (photograph in \cite{NORD3}):%
\begin{equation}
M_{6}=\left[ 
\begin{array}{cccccc}
28 & 4 & 3 & 31 & 35 & 10\vspace{-0.01in} \\ 
36 & 18 & 21 & 24 & 11 & 1\vspace{-0.01in} \\ 
7 & 23 & 12 & 17 & 22 & 30\vspace{-0.01in} \\ 
8 & 13 & 26 & 19 & 16 & 29\vspace{-0.01in} \\ 
5 & 20 & 15 & 14 & 25 & 32\vspace{-0.01in} \\ 
27 & 33 & 34 & 6 & 2 & 9\vspace{-0.01in}%
\end{array}%
\right]   \label{X6}
\end{equation}%
with auxiliary matrices from (\ref{RnQn}):%
\begin{equation}
Q_{6}=\left[ 
\begin{array}{cccccc}
4 & 0 & 0 & 5 & 5 & 1\vspace{-0.01in} \\ 
5 & 2 & 3 & 3 & 1 & 0\vspace{-0.01in} \\ 
1 & 3 & 1 & 2 & 3 & 4\vspace{-0.01in} \\ 
1 & 2 & 4 & 3 & 2 & 4\vspace{-0.01in} \\ 
0 & 3 & 2 & 2 & 4 & 5\vspace{-0.01in} \\ 
4 & 5 & 5 & 0 & 0 & 1\vspace{-0.01in}%
\end{array}%
\right] ,\quad R_{6}=\left[ 
\begin{array}{cccccc}
3 & 3 & 2 & 0 & 4 & 3\vspace{-0.01in} \\ 
5 & 5 & 2 & 5 & 4 & 0\vspace{-0.01in} \\ 
0 & 4 & 5 & 4 & 3 & 5\vspace{-0.01in} \\ 
1 & 0 & 1 & 0 & 3 & 4\vspace{-0.01in} \\ 
4 & 1 & 2 & 1 & 0 & 1\vspace{-0.01in} \\ 
2 & 2 & 3 & 5 & 1 & 2\vspace{-0.01in}%
\end{array}%
\right] .  \label{XQR6}
\end{equation}%
Again $Q_{6}$ and $R_{6}$ are magic with $\hat{m}=15$.

The method of constructing the above three $M_{6}$ squares is not apparent
and it does not appear to make use of auxiliary squares.\medskip

\noindent {\Large Order-8 Squares}\smallskip

We consider the Eulerian composition of three Franklin squares given by
Pasles \cite[pp. 236-237]{PASL}. Several of these squares also have been
decomposed by Morris \cite{MORR}. Here is Franklin's original order-8 square
published in 1769 \cite[p. 127]{PASL}:%
\begin{equation}
F_{8}=\left[ 
\begin{array}{cccccccc}
52 & 61 & 4 & 13 & 20 & 29 & 36 & 45\vspace{-0.01in} \\ 
14 & 3 & 62 & 51 & 46 & 35 & 30 & 19\vspace{-0.01in} \\ 
53 & 60 & 5 & 12 & 21 & 28 & 37 & 44\vspace{-0.01in} \\ 
11 & 6 & 59 & 54 & 43 & 38 & 27 & 22\vspace{-0.01in} \\ 
55 & 58 & 7 & 10 & 23 & 26 & 39 & 42\vspace{-0.01in} \\ 
9 & 8 & 57 & 56 & 41 & 40 & 25 & 24\vspace{-0.01in} \\ 
50 & 63 & 2 & 15 & 18 & 31 & 34 & 47\vspace{-0.01in} \\ 
16 & 1 & 64 & 49 & 48 & 33 & 32 & 17\vspace{-0.01in}%
\end{array}%
\right]   \label{F8}
\end{equation}%
with auxiliary matrices from (\ref{RnQn}) and Pasles \cite[p. 236]{PASL}: 
\begin{equation}
Q_{8}=\left[ 
\begin{array}{cccccccc}
\mathbf{6} & \mathbf{7} & \mathbf{0} & \mathbf{1} & \mathbf{2} & \mathbf{3}
& \mathbf{4} & \mathbf{5\vspace{-0.01in}} \\ 
1 & 0 & 7 & 6 & 5 & 4 & 3 & 2\vspace{-0.01in} \\ 
6 & 7 & 0 & 1 & 2 & 3 & 4 & 5\vspace{-0.01in} \\ 
1 & 0 & 7 & 6 & 5 & 4 & 3 & 2\vspace{-0.01in} \\ 
6 & 7 & 0 & 1 & 2 & 3 & 4 & 5\vspace{-0.01in} \\ 
1 & 0 & 7 & 6 & 5 & 4 & 3 & 2\vspace{-0.01in} \\ 
6 & 7 & 0 & 1 & 2 & 3 & 4 & 5\vspace{-0.01in} \\ 
1 & 0 & 7 & 6 & 5 & 4 & 3 & 2\vspace{-0.01in}%
\end{array}%
\right] ,\quad R_{8}=\left[ 
\begin{array}{cccccccc}
\mathbf{3} & 4 & 3 & 4 & 3 & 4 & 3 & 4\vspace{-0.01in} \\ 
5 & \mathbf{2} & 5 & 2 & 5 & 2 & 5 & 2\vspace{-0.01in} \\ 
4 & 3 & \mathbf{4} & 3 & 4 & 3 & 4 & 3\vspace{-0.01in} \\ 
2 & 5 & 2 & \mathbf{5} & 2 & 5 & 2 & 5\vspace{-0.01in} \\ 
6 & 1 & 6 & 1 & \mathbf{6} & 1 & 6 & 1\vspace{-0.01in} \\ 
0 & 7 & 0 & 7 & 0 & \mathbf{7} & 0 & 7\vspace{-0.01in} \\ 
1 & 6 & 1 & 6 & 1 & 6 & \mathbf{1} & 6\vspace{-0.01in} \\ 
7 & 0 & 7 & 0 & 7 & 0 & 7 & \mathbf{0}\vspace{-0.01in}%
\end{array}%
\right] .  \label{QR8}
\end{equation}%
These auxiliary matrices satisfy Franklin's three sum conditions with $%
\tilde{m}=28$. They follow from formulas given by Nordgren \cite{NORD1}
based on the form of the first row of $Q_{8}$ and the main diagonal of $%
R_{8}.$ Also, this $F_{8}$ can be constructed by the direct method of Jacobs 
\cite{JACO}.

Here is Franklin's order-8 pandiagonal magic square given by Pasles \cite[p.
207]{PASL}:%
\begin{equation}
F_{8}=\left[ 
\begin{array}{cccccccc}
2 & 57 & 6 & 61 & 8 & 63 & 4 & 59\vspace{-0.01in} \\ 
7 & 64 & 3 & 60 & 1 & 58 & 5 & 62\vspace{-0.01in} \\ 
49 & 10 & 53 & 14 & 55 & 16 & 51 & 12\vspace{-0.01in} \\ 
56 & 15 & 52 & 11 & 50 & 9 & 54 & 13\vspace{-0.01in} \\ 
42 & 17 & 46 & 21 & 48 & 23 & 44 & 19\vspace{-0.01in} \\ 
47 & 24 & 43 & 20 & 41 & 18 & 45 & 22\vspace{-0.01in} \\ 
25 & 34 & 29 & 38 & 31 & 40 & 27 & 36\vspace{-0.01in} \\ 
32 & 39 & 28 & 35 & 26 & 33 & 30 & 37\vspace{-0.01in}%
\end{array}%
\right] .  \label{FP8}
\end{equation}%
It satisfies Franklin's bent diagonal sum conditions and his $2\times 2$
subsquare sum conditions but not his half-row/column sum conditions. Its
auxiliary matrices from (\ref{RnQn}) and Pasles \cite[p. 237]{PASL} are:%
\begin{equation}
Q_{8}=\left[ 
\begin{array}{cccccccc}
0 & 7 & 0 & 7 & 0 & 7 & 0 & 7\vspace{-0.01in} \\ 
0 & 7 & 0 & 7 & 0 & 7 & 0 & 7\vspace{-0.01in} \\ 
6 & 1 & 6 & 1 & 6 & 1 & 6 & 1\vspace{-0.01in} \\ 
6 & 1 & 6 & 1 & 6 & 1 & 6 & 1\vspace{-0.01in} \\ 
5 & 2 & 5 & 2 & 5 & 2 & 5 & 2\vspace{-0.01in} \\ 
5 & 2 & 5 & 2 & 5 & 2 & 5 & 2\vspace{-0.01in} \\ 
3 & 4 & 3 & 4 & 3 & 4 & 3 & 4\vspace{-0.01in} \\ 
3 & 4 & 3 & 4 & 3 & 4 & 3 & 4\vspace{-0.01in}%
\end{array}%
\right] ,\quad R_{8}=\left[ 
\begin{array}{cccccccc}
1 & 0 & 5 & 4 & 7 & 6 & 3 & 2\vspace{-0.01in} \\ 
6 & 7 & 2 & 3 & 0 & 1 & 4 & 5\vspace{-0.01in} \\ 
0 & 1 & 4 & 5 & 6 & 7 & 2 & 3\vspace{-0.01in} \\ 
7 & 6 & 3 & 2 & 1 & 0 & 5 & 4\vspace{-0.01in} \\ 
1 & 0 & 5 & 4 & 7 & 6 & 3 & 2\vspace{-0.01in} \\ 
6 & 7 & 2 & 3 & 0 & 1 & 4 & 5\vspace{-0.01in} \\ 
0 & 1 & 4 & 5 & 6 & 7 & 2 & 3\vspace{-0.01in} \\ 
7 & 6 & 3 & 2 & 1 & 0 & 5 & 4\vspace{-0.01in}%
\end{array}%
\right] .  \label{QRP8}
\end{equation}%
These auxiliary matrices also satisfy two of Franklin's three sum conditions
with $\tilde{m}=28$ and they are pandiagonal.

Here is another of Franklin's order-8 Franklin square given by Pasles \cite[%
p. 169]{PASL}:%
\begin{equation}
F_{8}=\left[ 
\begin{array}{cccccccc}
17 & 47 & 30 & 36 & 21 & 43 & 26 & 40\vspace{-0.01in} \\ 
32 & 34 & 19 & 45 & 28 & 38 & 23 & 41\vspace{-0.01in} \\ 
33 & 31 & 46 & 20 & 37 & 27 & 42 & 24\vspace{-0.01in} \\ 
48 & 18 & 35 & 29 & 44 & 22 & 39 & 25\vspace{-0.01in} \\ 
49 & 15 & 62 & 4 & 53 & 11 & 58 & 8\vspace{-0.01in} \\ 
64 & 2 & 51 & 13 & 60 & 6 & 55 & 9\vspace{-0.01in} \\ 
1 & 63 & 14 & 52 & 5 & 59 & 10 & 56\vspace{-0.01in} \\ 
16 & 50 & 3 & 61 & 12 & 54 & 7 & 57\vspace{-0.01in}%
\end{array}%
\right]   \label{FP88}
\end{equation}%
with auxiliary matrices from (\ref{RnQn}) and Pasles \cite[p. 237]{PASL}:

\begin{equation}
Q_{8}=\left[ 
\begin{array}{cccccccc}
2 & 5 & 3 & 4 & 2 & 5 & 3 & 4\vspace{-0.01in} \\ 
3 & 4 & 2 & 5 & 3 & 4 & 2 & 5\vspace{-0.01in} \\ 
4 & 3 & 5 & 2 & 4 & 3 & 5 & 2\vspace{-0.01in} \\ 
5 & 2 & 4 & 3 & 5 & 2 & 4 & 3\vspace{-0.01in} \\ 
6 & 1 & 7 & 0 & 6 & 1 & 7 & 0\vspace{-0.01in} \\ 
7 & 0 & 6 & 1 & 7 & 0 & 6 & 1\vspace{-0.01in} \\ 
0 & 7 & 1 & 6 & 0 & 7 & 1 & 6\vspace{-0.01in} \\ 
1 & 6 & 0 & 7 & 1 & 6 & 0 & 7\vspace{-0.01in}%
\end{array}%
\right] ,\quad R_{8}=\left[ 
\begin{array}{cccccccc}
0 & 6 & 5 & 3 & 4 & 2 & 1 & 7\vspace{-0.01in} \\ 
7 & 1 & 2 & 4 & 3 & 5 & 6 & 0\vspace{-0.01in} \\ 
0 & 6 & 5 & 3 & 4 & 2 & 1 & 7\vspace{-0.01in} \\ 
7 & 1 & 2 & 4 & 3 & 5 & 6 & 0\vspace{-0.01in} \\ 
0 & 6 & 5 & 3 & 4 & 2 & 1 & 7\vspace{-0.01in} \\ 
7 & 1 & 2 & 4 & 3 & 5 & 6 & 0\vspace{-0.01in} \\ 
0 & 6 & 5 & 3 & 4 & 2 & 1 & 7\vspace{-0.01in} \\ 
7 & 1 & 2 & 4 & 3 & 5 & 6 & 0\vspace{-0.01in}%
\end{array}%
\right] .  \label{QRP88}
\end{equation}%
These auxiliary matrices also satisfy Franklin's three sum conditions with $%
\tilde{m}=28.$ This $Q_{8}$ is pandiagonal but $R_{8}$ is not and so $F_{8}$
of (\ref{FP88}) is not pandiagonal.

The elements of all of the above auxiliary matrices follow a distinctive
pattern which may be how Franklin constructed them. A direct construction of
the two Franklin matrices (\ref{FP8},\ref{FP88}) appears to be rather
difficult. The reader is invited to try!

The following order-8 pandiagonal Franklin magic square is given by
Schindel, et.al. \negthinspace \cite[Supplement \#2574]{SCH}: 
\begin{equation}
F_{8}=\left[ 
\begin{array}{cccccccc}
6 & 55 & 12 & 57 & 2 & 51 & 16 & 61\vspace{-0.01in} \\ 
11 & 58 & 5 & 56 & 15 & 62 & 1 & 52\vspace{-0.01in} \\ 
53 & 8 & 59 & 10 & 49 & 4 & 63 & 14\vspace{-0.01in} \\ 
60 & 9 & 54 & 7 & 64 & 13 & 50 & 3\vspace{-0.01in} \\ 
21 & 40 & 27 & 42 & 17 & 36 & 31 & 46\vspace{-0.01in} \\ 
28 & 41 & 22 & 39 & 32 & 45 & 18 & 35\vspace{-0.01in} \\ 
38 & 23 & 44 & 25 & 34 & 19 & 48 & 29\vspace{-0.01in} \\ 
43 & 26 & 37 & 24 & 47 & 30 & 33 & 20\vspace{-0.01in}%
\end{array}%
\right]   \label{F888}
\end{equation}%
$\allowbreak $with auxiliary matrices from (\ref{RnQn}):%
\begin{equation}
Q_{8}=\left[ 
\begin{array}{cccccccc}
0 & 6 & 1 & 7 & 0 & 6 & 1 & 7\vspace{-0.01in} \\ 
1 & 7 & 0 & 6 & 1 & 7 & 0 & 6\vspace{-0.01in} \\ 
6 & 0 & 7 & 1 & 6 & 0 & 7 & 1\vspace{-0.01in} \\ 
7 & 1 & 6 & 0 & 7 & 1 & 6 & 0\vspace{-0.01in} \\ 
2 & 4 & 3 & 5 & 2 & 4 & 3 & 5\vspace{-0.01in} \\ 
3 & 5 & 2 & 4 & 3 & 5 & 2 & 4\vspace{-0.01in} \\ 
4 & 2 & 5 & 3 & 4 & 2 & 5 & 3\vspace{-0.01in} \\ 
5 & 3 & 4 & 2 & 5 & 3 & 4 & 2\vspace{-0.01in}%
\end{array}%
\right] ,\quad R_{8}=\left[ 
\begin{array}{cccccccc}
5 & 6 & 3 & 0 & 1 & 2 & 7 & 4\vspace{-0.01in} \\ 
2 & 1 & 4 & 7 & 6 & 5 & 0 & 3\vspace{-0.01in} \\ 
4 & 7 & 2 & 1 & 0 & 3 & 6 & 5\vspace{-0.01in} \\ 
3 & 0 & 5 & 6 & 7 & 4 & 1 & 2\vspace{-0.01in} \\ 
4 & 7 & 2 & 1 & 0 & 3 & 6 & 5\vspace{-0.01in} \\ 
3 & 0 & 5 & 6 & 7 & 4 & 1 & 2\vspace{-0.01in} \\ 
5 & 6 & 3 & 0 & 1 & 2 & 7 & 4\vspace{-0.01in} \\ 
2 & 1 & 4 & 7 & 6 & 5 & 0 & 3\vspace{-0.01in}%
\end{array}%
\right] .  \label{QR888}
\end{equation}%
These auxiliary matrices also satisfy Franklin's three sum conditions with $%
\tilde{m}=28$ and they are pandiagonal. Although there is a pattern to their
elements, it is not clear how they could be systematically constructed.
Auxiliary matrices of other Franklin matrices from \cite{SCH} show similar
patterns.\newpage 

\noindent {\Large Order-16 Squares}\smallskip

Here is Franklin's original order-16 square published in1769 according to
Pasles \cite[p. 135]{PASL}:%
\begin{equation}
F_{16}=\!\left[ \!%
\begin{array}{cccccccccccccccc}
200 & 217 & 232 & 249 & 8 & 25 & 40 & 57 & 72 & 89 & 104 & 121 & 136 & 153 & 
168 & 185 \\ 
58 & 39 & 26 & 7 & 250 & 231 & 218 & 199 & 186 & 167 & 154 & 135 & 122 & 103
& 90 & 71 \\ 
198 & 219 & 230 & 251 & 6 & 27 & 38 & 59 & 70 & 91 & 102 & 123 & 134 & 155 & 
166 & 187 \\ 
60 & 37 & 28 & 5 & 252 & 229 & 220 & 197 & 188 & 165 & 156 & 133 & 124 & 101
& 92 & 69 \\ 
201 & 216 & 233 & 248 & 9 & 24 & 41 & 56 & 73 & 88 & 105 & 120 & 137 & 152 & 
169 & 184 \\ 
55 & 42 & 23 & 10 & 247 & 234 & 215 & 202 & 183 & 170 & 151 & 138 & 119 & 106
& 87 & 74 \\ 
203 & 214 & 235 & 246 & 11 & 22 & 43 & 54 & 75 & 86 & 107 & 118 & 139 & 150
& 171 & 182 \\ 
53 & 44 & 21 & 12 & 245 & 236 & 213 & 204 & 181 & 172 & 149 & 140 & 117 & 108
& 85 & 76 \\ 
205 & 212 & 237 & 244 & 13 & 20 & 45 & 52 & 77 & 84 & 109 & 116 & 141 & 148
& 173 & 180 \\ 
51 & 46 & 19 & 14 & 243 & 238 & 211 & 206 & 179 & 174 & 147 & 142 & 115 & 110
& 83 & 78 \\ 
207 & 210 & 239 & 242 & 15 & 18 & 47 & 50 & 79 & 82 & 111 & 114 & 143 & 146
& 175 & 178 \\ 
49 & 48 & 17 & 16 & 241 & 240 & 209 & 208 & 177 & 176 & 145 & 144 & 113 & 112
& 81 & 80 \\ 
196 & 221 & 228 & 253 & 4 & 29 & 36 & 61 & 68 & 93 & 100 & 125 & 132 & 157 & 
164 & 189 \\ 
62 & 35 & 30 & 3 & 254 & 227 & 222 & 195 & 190 & 163 & 158 & 131 & 126 & 99
& 94 & 67 \\ 
194 & 223 & 226 & 255 & 2 & 31 & 34 & 63 & 66 & 95 & 98 & 127 & 130 & 159 & 
162 & 191 \\ 
64 & 33 & 32 & 1 & 256 & 225 & 224 & 193 & 192 & 161 & 160 & 129 & 128 & 97
& 96 & 65%
\end{array}%
\!\right]   \label{F16}
\end{equation}%
with auxiliary matrices from (\ref{RnQn}) and \cite{MORR}:%
\begin{equation}
Q_{16}=\left[ 
\begin{array}{cccccccccccccccc}
12 & 13 & 14 & 15 & 0 & 1 & 2 & 3 & 4 & 5 & 6 & 7 & 8 & 9 & 10 & 11 \\ 
3 & 2 & 1 & 0 & 15 & 14 & 13 & 12 & 11 & 10 & 9 & 8 & 7 & 6 & 5 & 4 \\ 
12 & 13 & 14 & 15 & 0 & 1 & 2 & 3 & 4 & 5 & 6 & 7 & 8 & 9 & 10 & 11 \\ 
3 & 2 & 1 & 0 & 15 & 14 & 13 & 12 & 11 & 10 & 9 & 8 & 7 & 6 & 5 & 4 \\ 
\vdots  & \vdots  & \vdots  & \vdots  & \vdots  & \vdots  & \vdots  & \vdots 
& \vdots  & \vdots  & \vdots  & \vdots  & \vdots  & \vdots  & \vdots  & 
\vdots  \\ 
12 & 13 & 14 & 15 & 0 & 1 & 2 & 3 & 4 & 5 & 6 & 7 & 8 & 9 & 10 & 11 \\ 
3 & 2 & 1 & 0 & 15 & 14 & 13 & 12 & 11 & 10 & 9 & 8 & 7 & 6 & 5 & 4%
\end{array}%
\right] ,  \label{Q16}
\end{equation}%
\begin{equation}
R_{16}=\left[ 
\begin{array}{cccccccccccccccc}
7 & 8 & 7 & 8 & 7 & 8 & 7 & 8 & 7 & 8 & 7 & 8 & 7 & 8 & 7 & 8 \\ 
9 & 6 & 9 & 6 & 9 & 6 & 9 & 6 & 9 & 6 & 9 & 6 & 9 & 6 & 9 & 6 \\ 
5 & 10 & 5 & 10 & 5 & 10 & 5 & 10 & 5 & 10 & 5 & 10 & 5 & 10 & 5 & 10 \\ 
11 & 4 & 11 & 4 & 11 & 4 & 11 & 4 & 11 & 4 & 11 & 4 & 11 & 4 & 11 & 4 \\ 
8 & 7 & 8 & 7 & 8 & 7 & 8 & 7 & 8 & 7 & 8 & 7 & 8 & 7 & 8 & 7 \\ 
6 & 9 & 6 & 9 & 6 & 9 & 6 & 9 & 6 & 9 & 6 & 9 & 6 & 9 & 6 & 9 \\ 
10 & 5 & 10 & 5 & 10 & 5 & 10 & 5 & 10 & 5 & 10 & 5 & 10 & 5 & 10 & 5 \\ 
4 & 11 & 4 & 11 & 4 & 11 & 4 & 11 & 4 & 11 & 4 & 11 & 4 & 11 & 4 & 11 \\ 
12 & 3 & 12 & 3 & 12 & 3 & 12 & 3 & 12 & 3 & 12 & 3 & 12 & 3 & 12 & 3 \\ 
2 & 13 & 2 & 13 & 2 & 13 & 2 & 13 & 2 & 13 & 2 & 13 & 2 & 13 & 2 & 13 \\ 
14 & 1 & 14 & 1 & 14 & 1 & 14 & 1 & 14 & 1 & 14 & 1 & 14 & 1 & 14 & 1 \\ 
0 & 15 & 0 & 15 & 0 & 15 & 0 & 15 & 0 & 15 & 0 & 15 & 0 & 15 & 0 & 15 \\ 
3 & 12 & 3 & 12 & 3 & 12 & 3 & 12 & 3 & 12 & 3 & 12 & 3 & 12 & 3 & 12 \\ 
13 & 2 & 13 & 2 & 13 & 2 & 13 & 2 & 13 & 2 & 13 & 2 & 13 & 2 & 13 & 2 \\ 
1 & 14 & 1 & 14 & 1 & 14 & 1 & 14 & 1 & 14 & 1 & 14 & 1 & 14 & 1 & 14 \\ 
15 & 0 & 15 & 0 & 15 & 0 & 15 & 0 & 15 & 0 & 15 & 0 & 15 & 0 & 15 & 0%
\end{array}%
\right] .  \label{RR16}
\end{equation}%
These auxiliary matrices satisfy Franklin's three sum conditions with $%
\tilde{m}=120$. They can be constructed by Nordgren's formulas \cite{NORD1}.
The Franklin square $F_{16}$ also is constructed by Jacobs' direct method 
\cite{JACO}.

Here is Franklin's pandiagonal square according to Pasles \cite[p. 202]{PASL}%
:%
\begin{equation}
F_{16}=\left[ \!%
\begin{array}{cccccccccccccccc}
16 & 255 & 2 & 241 & 14 & 253 & 4 & 243 & 12 & 251 & 6 & 245 & 10 & 249 & 8
& 247 \\ 
1 & 242 & 15 & 256 & 3 & 244 & 13 & 254 & 5 & 246 & 11 & 252 & 7 & 248 & 9 & 
250 \\ 
240 & 31 & 226 & 17 & 238 & 29 & 228 & 19 & 236 & 27 & 230 & 21 & 234 & 25 & 
232 & 23 \\ 
225 & 18 & 239 & 32 & 227 & 20 & 237 & 30 & 229 & 22 & 235 & 28 & 231 & 24 & 
233 & 26 \\ 
223 & 48 & 209 & 34 & 221 & 46 & 211 & 36 & 219 & 44 & 213 & 38 & 217 & 42 & 
215 & 40 \\ 
210 & 33 & 224 & 47 & 212 & 35 & 222 & 45 & 214 & 37 & 220 & 43 & 216 & 39 & 
218 & 41 \\ 
63 & 208 & 49 & 194 & 61 & 206 & 51 & 196 & 59 & 204 & 53 & 198 & 57 & 202 & 
55 & 200 \\ 
50 & 193 & 64 & 207 & 52 & 195 & 62 & 205 & 54 & 197 & 60 & 203 & 56 & 199 & 
58 & 201 \\ 
80 & 191 & 66 & 177 & 78 & 189 & 68 & 179 & 76 & 187 & 70 & 181 & 74 & 185 & 
72 & 183 \\ 
65 & 178 & 79 & 192 & 67 & 180 & 77 & 190 & 69 & 182 & 75 & 188 & 71 & 184 & 
73 & 186 \\ 
176 & 95 & 162 & 81 & 174 & 93 & 164 & 83 & 172 & 91 & 166 & 85 & 170 & 89 & 
168 & 87 \\ 
161 & 82 & 175 & 96 & 163 & 84 & 173 & 94 & 165 & 86 & 171 & 92 & 167 & 88 & 
169 & 90 \\ 
159 & 112 & 145 & 98 & 157 & 110 & 147 & 100 & 155 & 108 & 149 & 102 & 153 & 
106 & 151 & 104 \\ 
146 & 97 & 160 & 111 & 148 & 99 & 158 & 109 & 150 & 101 & 156 & 107 & 152 & 
103 & 154 & 105 \\ 
127 & 144 & 113 & 130 & 125 & 142 & 115 & 132 & 123 & 140 & 117 & 134 & 121
& 138 & 119 & 136 \\ 
114 & 129 & 128 & 143 & 116 & 131 & 126 & 141 & 118 & 133 & 124 & 139 & 120
& 135 & 122 & 137%
\end{array}%
\!\right]   \label{F16p}
\end{equation}%
with auxiliary matrices from (\ref{RnQn}) and \cite{MORR}:%
\begin{equation}
Q_{16}=\left[ 
\begin{array}{cccccccccccccccc}
0 & 15 & 0 & 15 & 0 & 15 & 0 & 15 & 0 & 15 & 0 & 15 & 0 & 15 & 0 & 15\vspace{%
-0.01in} \\ 
0 & 15 & 0 & 15 & 0 & 15 & 0 & 15 & 0 & 15 & 0 & 15 & 0 & 15 & 0 & 15\vspace{%
-0.01in} \\ 
14 & 1 & 14 & 1 & 14 & 1 & 14 & 1 & 14 & 1 & 14 & 1 & 14 & 1 & 14 & 1\vspace{%
-0.01in} \\ 
14 & 1 & 14 & 1 & 14 & 1 & 14 & 1 & 14 & 1 & 14 & 1 & 14 & 1 & 14 & 1\vspace{%
-0.01in} \\ 
13 & 2 & 13 & 2 & 13 & 2 & 13 & 2 & 13 & 2 & 13 & 2 & 13 & 2 & 13 & 2\vspace{%
-0.01in} \\ 
13 & 2 & 13 & 2 & 13 & 2 & 13 & 2 & 13 & 2 & 13 & 2 & 13 & 2 & 13 & 2\vspace{%
-0.01in} \\ 
3 & 12 & 3 & 12 & 3 & 12 & 3 & 12 & 3 & 12 & 3 & 12 & 3 & 12 & 3 & 12\vspace{%
-0.01in} \\ 
3 & 12 & 3 & 12 & 3 & 12 & 3 & 12 & 3 & 12 & 3 & 12 & 3 & 12 & 3 & 12\vspace{%
-0.01in} \\ 
4 & 11 & 4 & 11 & 4 & 11 & 4 & 11 & 4 & 11 & 4 & 11 & 4 & 11 & 4 & 11\vspace{%
-0.01in} \\ 
4 & 11 & 4 & 11 & 4 & 11 & 4 & 11 & 4 & 11 & 4 & 11 & 4 & 11 & 4 & 11\vspace{%
-0.01in} \\ 
10 & 5 & 10 & 5 & 10 & 5 & 10 & 5 & 10 & 5 & 10 & 5 & 10 & 5 & 10 & 5\vspace{%
-0.01in} \\ 
10 & 5 & 10 & 5 & 10 & 5 & 10 & 5 & 10 & 5 & 10 & 5 & 10 & 5 & 10 & 5\vspace{%
-0.01in} \\ 
9 & 6 & 9 & 6 & 9 & 6 & 9 & 6 & 9 & 6 & 9 & 6 & 9 & 6 & 9 & 6\vspace{-0.01in}
\\ 
9 & 6 & 9 & 6 & 9 & 6 & 9 & 6 & 9 & 6 & 9 & 6 & 9 & 6 & 9 & 6\vspace{-0.01in}
\\ 
7 & 8 & 7 & 8 & 7 & 8 & 7 & 8 & 7 & 8 & 7 & 8 & 7 & 8 & 7 & 8\vspace{-0.01in}
\\ 
7 & 8 & 7 & 8 & 7 & 8 & 7 & 8 & 7 & 8 & 7 & 8 & 7 & 8 & 7 & 8\vspace{-0.01in}%
\end{array}%
\right] ,  \label{Q16p}
\end{equation}%
\begin{equation}
R_{16}=\left[ 
\begin{array}{cccccccccccccccc}
15 & 14 & 1 & 0 & 13 & 12 & 3 & 2 & 11 & 10 & 5 & 4 & 9 & 8 & 7 & 6\vspace{%
-0.01in} \\ 
0 & 1 & 14 & 15 & 2 & 3 & 12 & 13 & 4 & 5 & 10 & 11 & 6 & 7 & 8 & 9\vspace{%
-0.01in} \\ 
15 & 14 & 1 & 0 & 13 & 12 & 3 & 2 & 11 & 10 & 5 & 4 & 9 & 8 & 7 & 6\vspace{%
-0.01in} \\ 
0 & 1 & 14 & 15 & 2 & 3 & 12 & 13 & 4 & 5 & 10 & 11 & 6 & 7 & 8 & 9\vspace{%
-0.01in} \\ 
14 & 15 & 0 & 1 & 12 & 13 & 2 & 3 & 10 & 11 & 4 & 5 & 8 & 9 & 6 & 7\vspace{%
-0.01in} \\ 
1 & 0 & 15 & 14 & 3 & 2 & 13 & 12 & 5 & 4 & 11 & 10 & 7 & 6 & 9 & 8\vspace{%
-0.01in} \\ 
14 & 15 & 0 & 1 & 12 & 13 & 2 & 3 & 10 & 11 & 4 & 5 & 8 & 9 & 6 & 7\vspace{%
-0.01in} \\ 
1 & 0 & 15 & 14 & 3 & 2 & 13 & 12 & 5 & 4 & 11 & 10 & 7 & 6 & 9 & 8\vspace{%
-0.01in} \\ 
15 & 14 & 1 & 0 & 13 & 12 & 3 & 2 & 11 & 10 & 5 & 4 & 9 & 8 & 7 & 6\vspace{%
-0.01in} \\ 
0 & 1 & 14 & 15 & 2 & 3 & 12 & 13 & 4 & 5 & 10 & 11 & 6 & 7 & 8 & 9\vspace{%
-0.01in} \\ 
15 & 14 & 1 & 0 & 13 & 12 & 3 & 2 & 11 & 10 & 5 & 4 & 9 & 8 & 7 & 6\vspace{%
-0.01in} \\ 
0 & 1 & 14 & 15 & 2 & 3 & 12 & 13 & 4 & 5 & 10 & 11 & 6 & 7 & 8 & 9\vspace{%
-0.01in} \\ 
14 & 15 & 0 & 1 & 12 & 13 & 2 & 3 & 10 & 11 & 4 & 5 & 8 & 9 & 6 & 7\vspace{%
-0.01in} \\ 
1 & 0 & 15 & 14 & 3 & 2 & 13 & 12 & 5 & 4 & 11 & 10 & 7 & 6 & 9 & 8\vspace{%
-0.01in} \\ 
14 & 15 & 0 & 1 & 12 & 13 & 2 & 3 & 10 & 11 & 4 & 5 & 8 & 9 & 6 & 7\vspace{%
-0.01in} \\ 
1 & 0 & 15 & 14 & 3 & 2 & 13 & 12 & 5 & 4 & 11 & 10 & 7 & 6 & 9 & 8\vspace{%
-0.01in}%
\end{array}%
\right] .  \label{R16p}
\end{equation}%
These auxiliary matrices satisfy Franklin's three sum conditions with $%
\tilde{m}=120$ and they are pandiagonal. Again their elements follow a
distinct pattern.

Here is a new order-16 natural, pandiagonal, Franklin square constructed by
generalizing the order-8 auxiliary squares (\ref{QRP8}):%
\begin{equation}
F_{16}=\left[ \!%
\begin{array}{cccccccccccccccc}
2 & 241 & 6 & 245 & 10 & 249 & 14 & 253 & 16 & 255 & 12 & 251 & 8 & 247 & 4
& 243 \\ 
15 & 256 & 11 & 252 & 7 & 248 & 3 & 244 & 1 & 242 & 5 & 246 & 9 & 250 & 13 & 
254 \\ 
225 & 18 & 229 & 22 & 233 & 26 & 237 & 30 & 239 & 32 & 235 & 28 & 231 & 24 & 
227 & 20 \\ 
240 & 31 & 236 & 27 & 232 & 23 & 228 & 19 & 226 & 17 & 230 & 21 & 234 & 25 & 
238 & 29 \\ 
210 & 33 & 214 & 37 & 218 & 41 & 222 & 45 & 224 & 47 & 220 & 43 & 216 & 39 & 
212 & 35 \\ 
223 & 48 & 219 & 44 & 215 & 40 & 211 & 36 & 209 & 34 & 213 & 38 & 217 & 42 & 
221 & 46 \\ 
49 & 194 & 53 & 198 & 57 & 202 & 61 & 206 & 63 & 208 & 59 & 204 & 55 & 200 & 
51 & 196 \\ 
64 & 207 & 60 & 203 & 56 & 199 & 52 & 195 & 50 & 193 & 54 & 197 & 58 & 201 & 
62 & 205 \\ 
66 & 177 & 70 & 181 & 74 & 185 & 78 & 189 & 80 & 191 & 76 & 187 & 72 & 183 & 
68 & 179 \\ 
79 & 192 & 75 & 188 & 71 & 184 & 67 & 180 & 65 & 178 & 69 & 182 & 73 & 186 & 
77 & 190 \\ 
161 & 82 & 165 & 86 & 169 & 90 & 173 & 94 & 175 & 96 & 171 & 92 & 167 & 88 & 
163 & 84 \\ 
176 & 95 & 172 & 91 & 168 & 87 & 164 & 83 & 162 & 81 & 166 & 85 & 170 & 89 & 
174 & 93 \\ 
146 & 97 & 150 & 101 & 154 & 105 & 158 & 109 & 160 & 111 & 156 & 107 & 152 & 
103 & 148 & 99 \\ 
159 & 112 & 155 & 108 & 151 & 104 & 147 & 100 & 145 & 98 & 149 & 102 & 153 & 
106 & 157 & 110 \\ 
113 & 130 & 117 & 134 & 121 & 138 & 125 & 142 & 127 & 144 & 123 & 140 & 119
& 136 & 115 & 132 \\ 
128 & 143 & 124 & 139 & 120 & 135 & 116 & 131 & 114 & 129 & 118 & 133 & 122
& 137 & 126 & 141%
\end{array}%
\!\right] .  \label{F16new}
\end{equation}%
Further generalizations should lead to higher order-$8k$ squares of this
same type.

Here is a new order-16 natural Franklin square constructed by generalizing
the order-8 auxiliary squares (\ref{QRP88}):

\begin{equation}
F_{16}=\left[ \!%
\begin{array}{cccccccccccccccc}
65 & 191 & 94 & 164 & 69 & 187 & 90 & 168 & 73 & 183 & 86 & 172 & 77 & 179 & 
82 & 176 \\ 
96 & 162 & 67 & 189 & 92 & 166 & 71 & 185 & 88 & 170 & 75 & 181 & 84 & 174 & 
79 & 177 \\ 
97 & 159 & 126 & 132 & 101 & 155 & 122 & 136 & 105 & 151 & 118 & 140 & 109 & 
147 & 114 & 144 \\ 
128 & 130 & 99 & 157 & 124 & 134 & 103 & 153 & 120 & 138 & 107 & 149 & 116 & 
142 & 111 & 145 \\ 
129 & 127 & 158 & 100 & 133 & 123 & 154 & 104 & 137 & 119 & 150 & 108 & 141
& 115 & 146 & 112 \\ 
160 & 98 & 131 & 125 & 156 & 102 & 135 & 121 & 152 & 106 & 139 & 117 & 148 & 
110 & 143 & 113 \\ 
161 & 95 & 190 & 68 & 165 & 91 & 186 & 72 & 169 & 87 & 182 & 76 & 173 & 83 & 
178 & 80 \\ 
192 & 66 & 163 & 93 & 188 & 70 & 167 & 89 & 184 & 74 & 171 & 85 & 180 & 78 & 
175 & 81 \\ 
193 & 63 & 222 & 36 & 197 & 59 & 218 & 40 & 201 & 55 & 214 & 44 & 205 & 51 & 
210 & 48 \\ 
224 & 34 & 195 & 61 & 220 & 38 & 199 & 57 & 216 & 42 & 203 & 53 & 212 & 46 & 
207 & 49 \\ 
225 & 31 & 254 & 4 & 229 & 27 & 250 & 8 & 233 & 23 & 246 & 12 & 237 & 19 & 
242 & 16 \\ 
256 & 2 & 227 & 29 & 252 & 6 & 231 & 25 & 248 & 10 & 235 & 21 & 244 & 14 & 
239 & 17 \\ 
1 & 255 & 30 & 228 & 5 & 251 & 26 & 232 & 9 & 247 & 22 & 236 & 13 & 243 & 18
& 240 \\ 
32 & 226 & 3 & 253 & 28 & 230 & 7 & 249 & 24 & 234 & 11 & 245 & 20 & 238 & 15
& 241 \\ 
33 & 223 & 62 & 196 & 37 & 219 & 58 & 200 & 41 & 215 & 54 & 204 & 45 & 211 & 
50 & 208 \\ 
64 & 194 & 35 & 221 & 60 & 198 & 39 & 217 & 56 & 202 & 43 & 213 & 52 & 206 & 
47 & 209%
\end{array}%
\!\right] .
\end{equation}%
Further generalizations should lead to higher order-$8k$ squares of this
same type. Several other squares of the Franklin type are constructed by
Morris \cite{MORR} using Eulerian composition. Also, Hurkens \cite{HURK}
constructs several Franklin squares of various orders by direct
methods.\medskip

\noindent {\Large Conclusion}\smallskip

In the author's opinion and that of others \cite{MORR,PASL}, the distinct
pattern of elements of the auxiliary matrices $Q_{n}$ and $R_{n}$ considered
here indicate that Franklin probably used the Eulerian composite method of
construction of his squares of orders 8 and 16. It is difficult to see how
many of his squares and others presented here and in \cite{MORR,PASL} could
be constructed by a direct method. However, the referees of Nordgren's
article \cite{NORD1} contend that Franklin used a direct method which shows
that the question of Franklin's construction method may never be resolved to
everyone's satisfaction.\medskip

\textbf{Acknowledgement \ \ }I thank Peter Loly and Donald Morris for
helpful discussions and Adam Rogers for sending me the 4320 basic order-8
Franklin squares from \cite{SCH}.

\begin{center}
\textbf{Figure 1}\quad Order-24 Franklin square\newpage

\begin{align*}
Q_{24}& =\left[ 
\begin{array}{cccccccccccc}
18 & 19 & 20 & 21 & 22 & 23 & 0 & 1 & 2 & 3 & 4 & 5 \\ 
5 & 4 & 3 & 2 & 1 & 0 & 23 & 22 & 21 & 20 & 19 & 18 \\ 
18 & 19 & 20 & 21 & 22 & 23 & 0 & 1 & 2 & 3 & 4 & 5 \\ 
5 & 4 & 3 & 2 & 1 & 0 & 23 & 22 & 21 & 20 & 19 & 18 \\ 
18 & 19 & 20 & 21 & 22 & 23 & 0 & 1 & 2 & 3 & 4 & 5 \\ 
5 & 4 & 3 & 2 & 1 & 0 & 23 & 22 & 21 & 20 & 19 & 18 \\ 
18 & 19 & 20 & 21 & 22 & 23 & 0 & 1 & 2 & 3 & 4 & 5 \\ 
5 & 4 & 3 & 2 & 1 & 0 & 23 & 22 & 21 & 20 & 19 & 18 \\ 
18 & 19 & 20 & 21 & 22 & 23 & 0 & 1 & 2 & 3 & 4 & 5 \\ 
5 & 4 & 3 & 2 & 1 & 0 & 23 & 22 & 21 & 20 & 19 & 18 \\ 
18 & 19 & 20 & 21 & 22 & 23 & 0 & 1 & 2 & 3 & 4 & 5 \\ 
5 & 4 & 3 & 2 & 1 & 0 & 23 & 22 & 21 & 20 & 19 & 18 \\ 
18 & 19 & 20 & 21 & 22 & 23 & 0 & 1 & 2 & 3 & 4 & 5 \\ 
5 & 4 & 3 & 2 & 1 & 0 & 23 & 22 & 21 & 20 & 19 & 18 \\ 
18 & 19 & 20 & 21 & 22 & 23 & 0 & 1 & 2 & 3 & 4 & 5 \\ 
5 & 4 & 3 & 2 & 1 & 0 & 23 & 22 & 21 & 20 & 19 & 18 \\ 
18 & 19 & 20 & 21 & 22 & 23 & 0 & 1 & 2 & 3 & 4 & 5 \\ 
5 & 4 & 3 & 2 & 1 & 0 & 23 & 22 & 21 & 20 & 19 & 18 \\ 
18 & 19 & 20 & 21 & 22 & 23 & 0 & 1 & 2 & 3 & 4 & 5 \\ 
5 & 4 & 3 & 2 & 1 & 0 & 23 & 22 & 21 & 20 & 19 & 18 \\ 
18 & 19 & 20 & 21 & 22 & 23 & 0 & 1 & 2 & 3 & 4 & 5 \\ 
5 & 4 & 3 & 2 & 1 & 0 & 23 & 22 & 21 & 20 & 19 & 18 \\ 
18 & 19 & 20 & 21 & 22 & 23 & 0 & 1 & 2 & 3 & 4 & 5 \\ 
5 & 4 & 3 & 2 & 1 & 0 & 23 & 22 & 21 & 20 & 19 & 18%
\end{array}%
\right. \\
& \left. 
\begin{array}{cccccccccccc}
6 & 7 & 8 & 9 & 10 & 11 & 12 & 13 & 14 & 15 & 16 & 17 \\ 
17 & 16 & 15 & 14 & 13 & 12 & 11 & 10 & 9 & 8 & 7 & 6 \\ 
6 & 7 & 8 & 9 & 10 & 11 & 12 & 13 & 14 & 15 & 16 & 17 \\ 
17 & 16 & 15 & 14 & 13 & 12 & 11 & 10 & 9 & 8 & 7 & 6 \\ 
6 & 7 & 8 & 9 & 10 & 11 & 12 & 13 & 14 & 15 & 16 & 17 \\ 
17 & 16 & 15 & 14 & 13 & 12 & 11 & 10 & 9 & 8 & 7 & 6 \\ 
6 & 7 & 8 & 9 & 10 & 11 & 12 & 13 & 14 & 15 & 16 & 17 \\ 
17 & 16 & 15 & 14 & 13 & 12 & 11 & 10 & 9 & 8 & 7 & 6 \\ 
6 & 7 & 8 & 9 & 10 & 11 & 12 & 13 & 14 & 15 & 16 & 17 \\ 
17 & 16 & 15 & 14 & 13 & 12 & 11 & 10 & 9 & 8 & 7 & 6 \\ 
6 & 7 & 8 & 9 & 10 & 11 & 12 & 13 & 14 & 15 & 16 & 17 \\ 
17 & 16 & 15 & 14 & 13 & 12 & 11 & 10 & 9 & 8 & 7 & 6 \\ 
6 & 7 & 8 & 9 & 10 & 11 & 12 & 13 & 14 & 15 & 16 & 17 \\ 
17 & 16 & 15 & 14 & 13 & 12 & 11 & 10 & 9 & 8 & 7 & 6 \\ 
6 & 7 & 8 & 9 & 10 & 11 & 12 & 13 & 14 & 15 & 16 & 17 \\ 
17 & 16 & 15 & 14 & 13 & 12 & 11 & 10 & 9 & 8 & 7 & 6 \\ 
6 & 7 & 8 & 9 & 10 & 11 & 12 & 13 & 14 & 15 & 16 & 17 \\ 
17 & 16 & 15 & 14 & 13 & 12 & 11 & 10 & 9 & 8 & 7 & 6 \\ 
6 & 7 & 8 & 9 & 10 & 11 & 12 & 13 & 14 & 15 & 16 & 17 \\ 
17 & 16 & 15 & 14 & 13 & 12 & 11 & 10 & 9 & 8 & 7 & 6 \\ 
6 & 7 & 8 & 9 & 10 & 11 & 12 & 13 & 14 & 15 & 16 & 17 \\ 
17 & 16 & 15 & 14 & 13 & 12 & 11 & 10 & 9 & 8 & 7 & 6 \\ 
6 & 7 & 8 & 9 & 10 & 11 & 12 & 13 & 14 & 15 & 16 & 17 \\ 
17 & 16 & 15 & 14 & 13 & 12 & 11 & 10 & 9 & 8 & 7 & 6%
\end{array}%
\right]
\end{align*}

\textbf{Figure 2}\quad Quotient matrix for $F_{24}$\newpage 
\begin{align*}
R_{24}& =\left[ 
\begin{array}{cccccccccccc}
11 & 12 & 11 & 12 & 11 & 12 & 11 & 12 & 11 & 12 & 11 & 12 \\ 
13 & 10 & 13 & 10 & 13 & 10 & 13 & 10 & 13 & 10 & 13 & 10 \\ 
9 & 14 & 9 & 14 & 9 & 14 & 9 & 14 & 9 & 14 & 9 & 14 \\ 
15 & 8 & 15 & 8 & 15 & 8 & 15 & 8 & 15 & 8 & 15 & 8 \\ 
7 & 16 & 7 & 16 & 7 & 16 & 7 & 16 & 7 & 16 & 7 & 16 \\ 
17 & 6 & 17 & 6 & 17 & 6 & 17 & 6 & 17 & 6 & 17 & 6 \\ 
12 & 11 & 12 & 11 & 12 & 11 & 12 & 11 & 12 & 11 & 12 & 11 \\ 
10 & 13 & 10 & 13 & 10 & 13 & 10 & 13 & 10 & 13 & 10 & 13 \\ 
14 & 9 & 14 & 9 & 14 & 9 & 14 & 9 & 14 & 9 & 14 & 9 \\ 
8 & 15 & 8 & 15 & 8 & 15 & 8 & 15 & 8 & 15 & 8 & 15 \\ 
16 & 7 & 16 & 7 & 16 & 7 & 16 & 7 & 16 & 7 & 16 & 7 \\ 
6 & 17 & 6 & 17 & 6 & 17 & 6 & 17 & 6 & 17 & 6 & 17 \\ 
18 & 5 & 18 & 5 & 18 & 5 & 18 & 5 & 18 & 5 & 18 & 5 \\ 
4 & 19 & 4 & 19 & 4 & 19 & 4 & 19 & 4 & 19 & 4 & 19 \\ 
20 & 3 & 20 & 3 & 20 & 3 & 20 & 3 & 20 & 3 & 20 & 3 \\ 
2 & 21 & 2 & 21 & 2 & 21 & 2 & 21 & 2 & 21 & 2 & 21 \\ 
22 & 1 & 22 & 1 & 22 & 1 & 22 & 1 & 22 & 1 & 22 & 1 \\ 
0 & 23 & 0 & 23 & 0 & 23 & 0 & 23 & 0 & 23 & 0 & 23 \\ 
5 & 18 & 5 & 18 & 5 & 18 & 5 & 18 & 5 & 18 & 5 & 18 \\ 
19 & 4 & 19 & 4 & 19 & 4 & 19 & 4 & 19 & 4 & 19 & 4 \\ 
3 & 20 & 3 & 20 & 3 & 20 & 3 & 20 & 3 & 20 & 3 & 20 \\ 
21 & 2 & 21 & 2 & 21 & 2 & 21 & 2 & 21 & 2 & 21 & 2 \\ 
1 & 22 & 1 & 22 & 1 & 22 & 1 & 22 & 1 & 22 & 1 & 22 \\ 
23 & 0 & 23 & 0 & 23 & 0 & 23 & 0 & 23 & 0 & 23 & 0%
\end{array}%
\right. \\
& \left. 
\begin{array}{cccccccccccc}
11 & 12 & 11 & 12 & 11 & 12 & 11 & 12 & 11 & 12 & 11 & 12 \\ 
13 & 10 & 13 & 10 & 13 & 10 & 13 & 10 & 13 & 10 & 13 & 10 \\ 
9 & 14 & 9 & 14 & 9 & 14 & 9 & 14 & 9 & 14 & 9 & 14 \\ 
15 & 8 & 15 & 8 & 15 & 8 & 15 & 8 & 15 & 8 & 15 & 8 \\ 
7 & 16 & 7 & 16 & 7 & 16 & 7 & 16 & 7 & 16 & 7 & 16 \\ 
17 & 6 & 17 & 6 & 17 & 6 & 17 & 6 & 17 & 6 & 17 & 6 \\ 
12 & 11 & 12 & 11 & 12 & 11 & 12 & 11 & 12 & 11 & 12 & 11 \\ 
10 & 13 & 10 & 13 & 10 & 13 & 10 & 13 & 10 & 13 & 10 & 13 \\ 
14 & 9 & 14 & 9 & 14 & 9 & 14 & 9 & 14 & 9 & 14 & 9 \\ 
8 & 15 & 8 & 15 & 8 & 15 & 8 & 15 & 8 & 15 & 8 & 15 \\ 
16 & 7 & 16 & 7 & 16 & 7 & 16 & 7 & 16 & 7 & 16 & 7 \\ 
6 & 17 & 6 & 17 & 6 & 17 & 6 & 17 & 6 & 17 & 6 & 17 \\ 
18 & 5 & 18 & 5 & 18 & 5 & 18 & 5 & 18 & 5 & 18 & 5 \\ 
4 & 19 & 4 & 19 & 4 & 19 & 4 & 19 & 4 & 19 & 4 & 19 \\ 
20 & 3 & 20 & 3 & 20 & 3 & 20 & 3 & 20 & 3 & 20 & 3 \\ 
2 & 21 & 2 & 21 & 2 & 21 & 2 & 21 & 2 & 21 & 2 & 21 \\ 
22 & 1 & 22 & 1 & 22 & 1 & 22 & 1 & 22 & 1 & 22 & 1 \\ 
0 & 23 & 0 & 23 & 0 & 23 & 0 & 23 & 0 & 23 & 0 & 23 \\ 
5 & 18 & 5 & 18 & 5 & 18 & 5 & 18 & 5 & 18 & 5 & 18 \\ 
19 & 4 & 19 & 4 & 19 & 4 & 19 & 4 & 19 & 4 & 19 & 4 \\ 
3 & 20 & 3 & 20 & 3 & 20 & 3 & 20 & 3 & 20 & 3 & 20 \\ 
21 & 2 & 21 & 2 & 21 & 2 & 21 & 2 & 21 & 2 & 21 & 2 \\ 
1 & 22 & 1 & 22 & 1 & 22 & 1 & 22 & 1 & 22 & 1 & 22 \\ 
23 & 0 & 23 & 0 & 23 & 0 & 23 & 0 & 23 & 0 & 23 & 0%
\end{array}%
\right]
\end{align*}

\textbf{Figure 3}\quad Remainder matrix for $F_{24}$\newpage

\begin{equation*}
\left[ 
\begin{array}{ccccccccccccc}
1220 & 1261 & 1300 & 1341 & 1380 & 1421 & 1460 & 1501 & 1540 & 1581 & 20 & 61
& 100 \\ 
382 & 339 & 302 & 259 & 222 & 179 & 142 & 99 & 62 & 19 & 1582 & 1539 & 1502
\\ 
1218 & 1263 & 1298 & 1343 & 1378 & 1423 & 1458 & 1503 & 1538 & 1583 & 18 & 63
& 98 \\ 
384 & 337 & 304 & 257 & 224 & 177 & 144 & 97 & 64 & 17 & 1584 & 1537 & 1504
\\ 
1216 & 1265 & 1296 & 1345 & 1376 & 1425 & 1456 & 1505 & 1536 & 1585 & 16 & 65
& 96 \\ 
386 & 335 & 306 & 255 & 226 & 175 & 146 & 95 & 66 & 15 & 1586 & 1535 & 1506
\\ 
1214 & 1267 & 1294 & 1347 & 1374 & 1427 & 1454 & 1507 & 1534 & 1587 & 14 & 67
& 94 \\ 
388 & 333 & 308 & 253 & 228 & 173 & 148 & 93 & 68 & 13 & 1588 & 1533 & 1508
\\ 
1212 & 1269 & 1292 & 1349 & 1372 & 1429 & 1452 & 1509 & 1532 & 1589 & 12 & 69
& 92 \\ 
390 & 331 & 310 & 251 & 230 & 171 & 150 & 91 & 70 & 11 & 1590 & 1531 & 1510
\\ 
1221 & 1260 & 1301 & 1340 & 1381 & 1420 & 1461 & 1500 & 1541 & 1580 & 21 & 60
& 101 \\ 
379 & 342 & 299 & 262 & 219 & 182 & 139 & 102 & 59 & 22 & 1579 & 1542 & 1499
\\ 
1223 & 1258 & 1303 & 1338 & 1383 & 1418 & 1463 & 1498 & 1543 & 1578 & 23 & 58
& 103 \\ 
377 & 344 & 297 & 264 & 217 & 184 & 137 & 104 & 57 & 24 & 1577 & 1544 & 1497
\\ 
1225 & 1256 & 1305 & 1336 & 1385 & 1416 & 1465 & 1496 & 1545 & 1576 & 25 & 56
& 105 \\ 
375 & 346 & 295 & 266 & 215 & 186 & 135 & 106 & 55 & 26 & 1575 & 1546 & 1495
\\ 
1227 & 1254 & 1307 & 1334 & 1387 & 1414 & 1467 & 1494 & 1547 & 1574 & 27 & 54
& 107 \\ 
373 & 348 & 293 & 268 & 213 & 188 & 133 & 108 & 53 & 28 & 1573 & 1548 & 1493
\\ 
1229 & 1252 & 1309 & 1332 & 1389 & 1412 & 1469 & 1492 & 1549 & 1572 & 29 & 52
& 109 \\ 
371 & 350 & 291 & 270 & 211 & 190 & 131 & 110 & 51 & 30 & 1571 & 1550 & 1491
\\ 
1231 & 1250 & 1311 & 1330 & 1391 & 1410 & 1471 & 1490 & 1551 & 1570 & 31 & 50
& 111 \\ 
369 & 352 & 289 & 272 & 209 & 192 & 129 & 112 & 49 & 32 & 1569 & 1552 & 1489
\\ 
1233 & 1248 & 1313 & 1328 & 1393 & 1408 & 1473 & 1488 & 1553 & 1568 & 33 & 48
& 113 \\ 
367 & 354 & 287 & 274 & 207 & 194 & 127 & 114 & 47 & 34 & 1567 & 1554 & 1487
\\ 
1235 & 1246 & 1315 & 1326 & 1395 & 1406 & 1475 & 1486 & 1555 & 1566 & 35 & 46
& 115 \\ 
365 & 356 & 285 & 276 & 205 & 196 & 125 & 116 & 45 & 36 & 1565 & 1556 & 1485
\\ 
1237 & 1244 & 1317 & 1324 & 1397 & 1404 & 1477 & 1484 & 1557 & 1564 & 37 & 44
& 117 \\ 
363 & 358 & 283 & 278 & 203 & 198 & 123 & 118 & 43 & 38 & 1563 & 1558 & 1483
\\ 
1239 & 1242 & 1319 & 1322 & 1399 & 1402 & 1479 & 1482 & 1559 & 1562 & 39 & 42
& 119 \\ 
361 & 360 & 281 & 280 & 201 & 200 & 121 & 120 & 41 & 40 & 1561 & 1560 & 1481
\\ 
1210 & 1271 & 1290 & 1351 & 1370 & 1431 & 1450 & 1511 & 1530 & 1591 & 10 & 71
& 90 \\ 
392 & 329 & 312 & 249 & 232 & 169 & 152 & 89 & 72 & 9 & 1592 & 1529 & 1512
\\ 
1208 & 1273 & 1288 & 1353 & 1368 & 1433 & 1448 & 1513 & 1528 & 1593 & 8 & 73
& 88 \\ 
394 & 327 & 314 & 247 & 234 & 167 & 154 & 87 & 74 & 7 & 1594 & 1527 & 1514
\\ 
1206 & 1275 & 1286 & 1355 & 1366 & 1435 & 1446 & 1515 & 1526 & 1595 & 6 & 75
& 86 \\ 
396 & 325 & 316 & 245 & 236 & 165 & 156 & 85 & 76 & 5 & 1596 & 1525 & 1516
\\ 
1204 & 1277 & 1284 & 1357 & 1364 & 1437 & 1444 & 1517 & 1524 & 1597 & 4 & 77
& 84 \\ 
398 & 323 & 318 & 243 & 238 & 163 & 158 & 83 & 78 & 3 & 1598 & 1523 & 1518
\\ 
1202 & 1279 & 1282 & 1359 & 1362 & 1439 & 1442 & 1519 & 1522 & 1599 & 2 & 79
& 82 \\ 
400 & 321 & 320 & 241 & 240 & 161 & 160 & 81 & 80 & 1 & 1600 & 1521 & 1520%
\end{array}%
\right. \smallskip
\end{equation*}

\textbf{Figure 4a}\quad Order-40 Franklin matrix - left third\newpage

\begin{equation*}
\begin{array}{ccccccccccccc}
141 & 180 & 221 & 260 & 301 & 340 & 381 & 420 & 461 & 500 & 541 & 580 & 621
\\ 
1459 & 1422 & 1379 & 1342 & 1299 & 1262 & 1219 & 1182 & 1139 & 1102 & 1059 & 
1022 & 979 \\ 
143 & 178 & 223 & 258 & 303 & 338 & 383 & 418 & 463 & 498 & 543 & 578 & 623
\\ 
1457 & 1424 & 1377 & 1344 & 1297 & 1264 & 1217 & 1184 & 1137 & 1104 & 1057 & 
1024 & 977 \\ 
145 & 176 & 225 & 256 & 305 & 336 & 385 & 416 & 465 & 496 & 545 & 576 & 625
\\ 
1455 & 1426 & 1375 & 1346 & 1295 & 1266 & 1215 & 1186 & 1135 & 1106 & 1055 & 
1026 & 975 \\ 
147 & 174 & 227 & 254 & 307 & 334 & 387 & 414 & 467 & 494 & 547 & 574 & 627
\\ 
1453 & 1428 & 1373 & 1348 & 1293 & 1268 & 1213 & 1188 & 1133 & 1108 & 1053 & 
1028 & 973 \\ 
149 & 172 & 229 & 252 & 309 & 332 & 389 & 412 & 469 & 492 & 549 & 572 & 629
\\ 
1451 & 1430 & 1371 & 1350 & 1291 & 1270 & 1211 & 1190 & 1131 & 1110 & 1051 & 
1030 & 971 \\ 
140 & 181 & 220 & 261 & 300 & 341 & 380 & 421 & 460 & 501 & 540 & 581 & 620
\\ 
1462 & 1419 & 1382 & 1339 & 1302 & 1259 & 1222 & 1179 & 1142 & 1099 & 1062 & 
1019 & 982 \\ 
138 & 183 & 218 & 263 & 298 & 343 & 378 & 423 & 458 & 503 & 538 & 583 & 618
\\ 
1464 & 1417 & 1384 & 1337 & 1304 & 1257 & 1224 & 1177 & 1144 & 1097 & 1064 & 
1017 & 984 \\ 
136 & 185 & 216 & 265 & 296 & 345 & 376 & 425 & 456 & 505 & 536 & 585 & 616
\\ 
1466 & 1415 & 1386 & 1335 & 1306 & 1255 & 1226 & 1175 & 1146 & 1095 & 1066 & 
1015 & 986 \\ 
134 & 187 & 214 & 267 & 294 & 347 & 374 & 427 & 454 & 507 & 534 & 587 & 614
\\ 
1468 & 1413 & 1388 & 1333 & 1308 & 1253 & 1228 & 1173 & 1148 & 1093 & 1068 & 
1013 & 988 \\ 
132 & 189 & 212 & 269 & 292 & 349 & 372 & 429 & 452 & 509 & 532 & 589 & 612
\\ 
1470 & 1411 & 1390 & 1331 & 1310 & 1251 & 1230 & 1171 & 1150 & 1091 & 1070 & 
1011 & 990 \\ 
130 & 191 & 210 & 271 & 290 & 351 & 370 & 431 & 450 & 511 & 530 & 591 & 610
\\ 
1472 & 1409 & 1392 & 1329 & 1312 & 1249 & 1232 & 1169 & 1152 & 1089 & 1072 & 
1009 & 992 \\ 
128 & 193 & 208 & 273 & 288 & 353 & 368 & 433 & 448 & 513 & 528 & 593 & 608
\\ 
1474 & 1407 & 1394 & 1327 & 1314 & 1247 & 1234 & 1167 & 1154 & 1087 & 1074 & 
1007 & 994 \\ 
126 & 195 & 206 & 275 & 286 & 355 & 366 & 435 & 446 & 515 & 526 & 595 & 606
\\ 
1476 & 1405 & 1396 & 1325 & 1316 & 1245 & 1236 & 1165 & 1156 & 1085 & 1076 & 
1005 & 996 \\ 
124 & 197 & 204 & 277 & 284 & 357 & 364 & 437 & 444 & 517 & 524 & 597 & 604
\\ 
1478 & 1403 & 1398 & 1323 & 1318 & 1243 & 1238 & 1163 & 1158 & 1083 & 1078 & 
1003 & 998 \\ 
122 & 199 & 202 & 279 & 282 & 359 & 362 & 439 & 442 & 519 & 522 & 599 & 602
\\ 
1480 & 1401 & 1400 & 1321 & 1320 & 1241 & 1240 & 1161 & 1160 & 1081 & 1080 & 
1001 & 1000 \\ 
151 & 170 & 231 & 250 & 311 & 330 & 391 & 410 & 471 & 490 & 551 & 570 & 631
\\ 
1449 & 1432 & 1369 & 1352 & 1289 & 1272 & 1209 & 1192 & 1129 & 1112 & 1049 & 
1032 & 969 \\ 
153 & 168 & 233 & 248 & 313 & 328 & 393 & 408 & 473 & 488 & 553 & 568 & 633
\\ 
1447 & 1434 & 1367 & 1354 & 1287 & 1274 & 1207 & 1194 & 1127 & 1114 & 1047 & 
1034 & 967 \\ 
155 & 166 & 235 & 246 & 315 & 326 & 395 & 406 & 475 & 486 & 555 & 566 & 635
\\ 
1445 & 1436 & 1365 & 1356 & 1285 & 1276 & 1205 & 1196 & 1125 & 1116 & 1045 & 
1036 & 965 \\ 
157 & 164 & 237 & 244 & 317 & 324 & 397 & 404 & 477 & 484 & 557 & 564 & 637
\\ 
1443 & 1438 & 1363 & 1358 & 1283 & 1278 & 1203 & 1198 & 1123 & 1118 & 1043 & 
1038 & 963 \\ 
159 & 162 & 239 & 242 & 319 & 322 & 399 & 402 & 479 & 482 & 559 & 562 & 639
\\ 
1441 & 1440 & 1361 & 1360 & 1281 & 1280 & 1201 & 1200 & 1121 & 1120 & 1041 & 
1040 & 961%
\end{array}%
\end{equation*}%
$\smallskip \smallskip $\textbf{Figure 4b}\quad Order-40 Franklin matrix -
middle third\newpage

\begin{equation*}
\left. 
\begin{array}{cccccccccccccc}
660 & 701 & 740 & 781 & 820 & 861 & 900 & 941 & 980 & 1021 & 1060 & 1101 & 
1140 & 1181 \\ 
942 & 899 & 862 & 819 & 782 & 739 & 702 & 659 & 622 & 579 & 542 & 499 & 462
& 419 \\ 
658 & 703 & 738 & 783 & 818 & 863 & 898 & 943 & 978 & 1023 & 1058 & 1103 & 
1138 & 1183 \\ 
944 & 897 & 864 & 817 & 784 & 737 & 704 & 657 & 624 & 577 & 544 & 497 & 464
& 417 \\ 
656 & 705 & 736 & 785 & 816 & 865 & 896 & 945 & 976 & 1025 & 1056 & 1105 & 
1136 & 1185 \\ 
946 & 895 & 866 & 815 & 786 & 735 & 706 & 655 & 626 & 575 & 546 & 495 & 466
& 415 \\ 
654 & 707 & 734 & 787 & 814 & 867 & 894 & 947 & 974 & 1027 & 1054 & 1107 & 
1134 & 1187 \\ 
948 & 893 & 868 & 813 & 788 & 733 & 708 & 653 & 628 & 573 & 548 & 493 & 468
& 413 \\ 
652 & 709 & 732 & 789 & 812 & 869 & 892 & 949 & 972 & 1029 & 1052 & 1109 & 
1132 & 1189 \\ 
950 & 891 & 870 & 811 & 790 & 731 & 710 & 651 & 630 & 571 & 550 & 491 & 470
& 411 \\ 
661 & 700 & 741 & 780 & 821 & 860 & 901 & 940 & 981 & 1020 & 1061 & 1100 & 
1141 & 1180 \\ 
939 & 902 & 859 & 822 & 779 & 742 & 699 & 662 & 619 & 582 & 539 & 502 & 459
& 422 \\ 
663 & 698 & 743 & 778 & 823 & 858 & 903 & 938 & 983 & 1018 & 1063 & 1098 & 
1143 & 1178 \\ 
937 & 904 & 857 & 824 & 777 & 744 & 697 & 664 & 617 & 584 & 537 & 504 & 457
& 424 \\ 
665 & 696 & 745 & 776 & 825 & 856 & 905 & 936 & 985 & 1016 & 1065 & 1096 & 
1145 & 1176 \\ 
935 & 906 & 855 & 826 & 775 & 746 & 695 & 666 & 615 & 586 & 535 & 506 & 455
& 426 \\ 
667 & 694 & 747 & 774 & 827 & 854 & 907 & 934 & 987 & 1014 & 1067 & 1094 & 
1147 & 1174 \\ 
933 & 908 & 853 & 828 & 773 & 748 & 693 & 668 & 613 & 588 & 533 & 508 & 453
& 428 \\ 
669 & 692 & 749 & 772 & 829 & 852 & 909 & 932 & 989 & 1012 & 1069 & 1092 & 
1149 & 1172 \\ 
931 & 910 & 851 & 830 & 771 & 750 & 691 & 670 & 611 & 590 & 531 & 510 & 451
& 430 \\ 
671 & 690 & 751 & 770 & 831 & 850 & 911 & 930 & 991 & 1010 & 1071 & 1090 & 
1151 & 1170 \\ 
929 & 912 & 849 & 832 & 769 & 752 & 689 & 672 & 609 & 592 & 529 & 512 & 449
& 432 \\ 
673 & 688 & 753 & 768 & 833 & 848 & 913 & 928 & 993 & 1008 & 1073 & 1088 & 
1153 & 1168 \\ 
927 & 914 & 847 & 834 & 767 & 754 & 687 & 674 & 607 & 594 & 527 & 514 & 447
& 434 \\ 
675 & 686 & 755 & 766 & 835 & 846 & 915 & 926 & 995 & 1006 & 1075 & 1086 & 
1155 & 1166 \\ 
925 & 916 & 845 & 836 & 765 & 756 & 685 & 676 & 605 & 596 & 525 & 516 & 445
& 436 \\ 
677 & 684 & 757 & 764 & 837 & 844 & 917 & 924 & 997 & 1004 & 1077 & 1084 & 
1157 & 1164 \\ 
923 & 918 & 843 & 838 & 763 & 758 & 683 & 678 & 603 & 598 & 523 & 518 & 443
& 438 \\ 
679 & 682 & 759 & 762 & 839 & 842 & 919 & 922 & 999 & 1002 & 1079 & 1082 & 
1159 & 1162 \\ 
921 & 920 & 841 & 840 & 761 & 760 & 681 & 680 & 601 & 600 & 521 & 520 & 441
& 440 \\ 
650 & 711 & 730 & 791 & 810 & 871 & 890 & 951 & 970 & 1031 & 1050 & 1111 & 
1130 & 1191 \\ 
952 & 889 & 872 & 809 & 792 & 729 & 712 & 649 & 632 & 569 & 552 & 489 & 472
& 409 \\ 
648 & 713 & 728 & 793 & 808 & 873 & 888 & 953 & 968 & 1033 & 1048 & 1113 & 
1128 & 1193 \\ 
954 & 887 & 874 & 807 & 794 & 727 & 714 & 647 & 634 & 567 & 554 & 487 & 474
& 407 \\ 
646 & 715 & 726 & 795 & 806 & 875 & 886 & 955 & 966 & 1035 & 1046 & 1115 & 
1126 & 1195 \\ 
956 & 885 & 876 & 805 & 796 & 725 & 716 & 645 & 636 & 565 & 556 & 485 & 476
& 405 \\ 
644 & 717 & 724 & 797 & 804 & 877 & 884 & 957 & 964 & 1037 & 1044 & 1117 & 
1124 & 1197 \\ 
958 & 883 & 878 & 803 & 798 & 723 & 718 & 643 & 638 & 563 & 558 & 483 & 478
& 403 \\ 
642 & 719 & 722 & 799 & 802 & 879 & 882 & 959 & 962 & 1039 & 1042 & 1119 & 
1122 & 1199 \\ 
960 & 881 & 880 & 801 & 800 & 721 & 720 & 641 & 640 & 561 & 560 & 481 & 480
& 401%
\end{array}%
\right] 
\end{equation*}%
$\smallskip \smallskip $\textbf{Figure 4c}\quad Order-40 Franklin matrix -
right third
\end{center}

\end{document}